# Generalization of the Hardy-Littlewood conjecture to almost-prime number tuples

Victor Volfson

ABSTRACT. The article presents a generalization of the classical Hardy-Littlewood conjecture concerning the density of prime tuples to the case of tuples consisting of almost-prime numbers (numbers with a specified quantity of prime divisors). The work investigates tuples of natural numbers where each element is subject to an individual factorization requirement. A proposed asymptotic formula for the quantity of such tuples is presented, where the density is determined by the product of two constants: the standard Selberg constant, which depends solely on the tuple pattern, and a correction factor, which depends only on the set of requirements for the number of prime divisors at each position in the tuple. The author proves that the admissibility of a pattern for prime numbers implies its admissibility for almost-prime numbers. The principle of symmetry is established - the correction factor depends only on the multiset of requirements, not on the order of elements within the tuple. An empirical-analytical method for calculating the correction factor is developed, based on the invariance of the Selberg constant under pattern "stretching." The method is tested on tuples of small length (pairs and triples), for which tables of calculated coefficients with high accuracy are provided. The method is justified in the work.

## 1. INTRODUCTION

It is known that any natural number $n$ can be represented as:

$$n = \prod_{i=1}^{m} p_i^{a_i}, \qquad (1.1)$$

where $p_i$ is the $i$-th prime number, and $a_i$ is a natural number. This decomposition of a natural number is called a factorization.

Hardy and Littlewood first formulated the principle in 1923 that the density of prime tuples is determined by the product, over prime numbers $p$, of the probabilities that all elements of the tuple are not divisible by $p$ [1]. This can be generalized to the case of natural numbers with some factorization.

Let us call the $k$-almost prime number a natural number for which:

$$k = \sum_{i=1}^{m} a_i. \qquad (1.2)$$

Chen in 1973 [2] studied the first non-trivial case of this generalization for tuples $(n, n+2)$ where the first number is prime ($k=1$) and the second is 2-almost prime ($k=2$) in accordance with (1.1) and (1.2).

Heath-Brown studied the distribution of various near-primes on short intervals and progressions in 1978 [3]. His methods show that non-trivial estimates are obtained for "sparse" tuples with different factorizations (1.1).

Halberstam & Richert developed a general sieve theory in 1975 [4], which allows one to work with weighted sums of the form:

$$\sum_{n} \prod_{i=1}^{m} f_i(n + h_i),$$

where $f_i$ are multiplicative functions that satisfy certain factorization conditions.

Green & Tao proved that pseudorandom weight functions (like the Weyl function for primes) preserve arithmetic progressions in 2008 [5]. This opens the way to proving the existence of tuples where different positions have different factorizations.

Tao & Ziegler showed that the structure theorems hold for polynomial configurations in 2008 [6]. This corresponds to considering tuples of the form $(P_1(n), ..., P_m(n))$ with a different factorization of each polynomial value.

Maynard showed that the optimal sieve method allows one to separately evaluate the contributions of different positions in a tuple in 2015 [7], which is critically important for cases of tuples with different factorizations of natural numbers.

The classical Hardy-Littlewood conjecture (on $m$ prime tuples) states [1] that if the number of tuples with a pattern $H$ is infinite (an admissible pattern), then the asymptotic number of such tuples, on the interval up to $x$ with $x \to \infty$, is:

$$\pi_H(x) \sim \mathfrak{S}(H) \frac{x}{(\log x)^m}, \qquad (1.3)$$

where $\mathfrak{S}(H)$ is Selberg's constant.

Let's consider a tuple of natural numbers, where each number has its own fixed factorization type (1.1). Such a tuple will correspond to an integer vector $K = (k_1, k_2, ..., k_m)$, which defines the factorization requirements for each natural number in this tuple.

We will consider a generalization for twin primes, i.e., a tuple of length 2 $(n, n+2)$ as an example 1:

- the classical Hardy-Littlewood conjecture considers the case where both numbers in the tuple are prime and the integer vector has the form: $K = (1,1)$;

- in Chen's theorem [2], the value $n$ is prime, and $n+2$ is 2-almost prime and the integer vector has the form: $K = (1,2)$;

- a new generalization: $n$ is $k_1$-almost prime, is $k_2$-almost prime and the integer vector has the form: $K = (k_1, k_2)$.

Let's consider an arithmetic progression, i.e. a tuple $(n, n+d, n+2d, ..., n+(m-1)d)$, as an example 2:

- all numbers are prime in the Green-Tao theorem [5], therefore the integer vector has the form: $K = (1,1,...,1)$;

- new generalization: all numbers are $k_i$ almost prime and an integer vector has the form: $K = (k_1, k_2, ..., k_m)$.

Now we examine the asymptotics and distribution of the number of distinct almost prime numbers.

Let us consider two cases when a natural number $n$ has exactly $k$ different prime divisors, i.e. $\omega(n) = k$ and has $k$ prime divisors taking into account the multiplicity, i.e. $\Omega(n) = k$.

General asymptotics were obtained by Norton in 1976 and Tenenbaum in 1995 [8].

The asymptotics for the first case are as follows:

$$N_k(x) = \frac{x}{\log x} \cdot \frac{(\log \log x)^{k-1}}{(k-1)!} \left[ 1 + \frac{(k-1)\gamma - \frac{k(k-1)}{2} + C_1}{\log \log x} + O_k\left(\frac{1}{(\log \log x)^2}\right) \right], \qquad (1.4)$$

where $C_1 = \sum_p \left( \log(1 - \frac{1}{p}) + \frac{1}{p} \right) \approx -0.315718$.

The asymptotics for the second case are as follows:

$$\pi_k(x) = \frac{x}{\log x} \cdot \frac{(\log \log x)^{k-1}}{(k-1)!} \left[ 1 + \frac{(k-1)\gamma - \frac{k(k-1)}{2} + C_2}{\log \log x} + O_k\left(\frac{1}{(\log \log x)^2}\right) \right]. \qquad (1.5)$$

where $C_2 = \sum_p \left( \log(1 - \frac{1}{p}) + \frac{1}{p-1} \right) \approx 0.754916$.

Based on formulas (1.4) and (1.5), we can conclude that the leading terms of the asymptotics for both cases are the same:

$$N_k(x) \sim \pi_k(x) \sim \frac{x}{\log x} \cdot \frac{(\log \log x)^{k-1}}{(k-1)!}. \qquad (1.6)$$

According to the asymptotic formula (1.6), the number of $r$- almost prime numbers $\pi_r(x)$ for a fixed $x$ does not necessarily increase with increasing $r$.

Formula:

$$\pi_k(x) \sim \frac{x}{\log x} \cdot \frac{(\log \log x)^{k-1}}{(k-1)!}$$

shows that the behavior of $\pi_k(x)$, as a function of $k$ (for a fixed $x$), has a maximum, since the ratio of successive values is:

$$\frac{\pi_{k+1}(x)}{\pi_k(x)} \sim \frac{\log \log x}{k}. \qquad (1.7)$$

Based on (1.7):

- if $k < \log \log x$, then the ratio is greater than 1 and $\pi_k(x)$ increases with increasing $k$;

- if $k > \log \log x$, then the ratio is less than 1 and $\pi_k(x)$ decreases with increasing $k$.

The maximum is reached at $k \approx \log \log x$, which is consistent with the Hardy-Ramanujan theorem on the average number of prime divisors of a number.

Therefore (for fixed $x$) the number of $k$-almost primes first increases with $k$, reaches a maximum at $k \approx \log \log x$, and then decreases. For example, we have $\log \log x \approx 5.4$ for $x = 10^{100}$, so the maximum will be around $k = 5$ or 6.

The distribution of the number of prime divisors of a number is asymptotically normal in the interval around $x$ with the mean $\mu \sim \log\log x$ and variance $\sigma^2 \sim \log\log x$ (Erdős-Kac theorem).

However, (for finite $x$) the distribution is more accurately described by the Poisson distribution with parameter $\lambda = \log\log x$. The Poisson distribution approaches the normal distribution for large $\lambda$, but the rate of convergence is slow.

The distribution will be close to normal if the parameter $\lambda = \log\log x$ is large enough.

It is generally assumed, that the Poisson distribution is already fairly well approximated by the normal distribution for $\lambda > 10$. However, it is required often $\lambda > 20$ or even more to achieve symmetry and eliminate asymmetry (skewness).

Let us find $x$ for which $\log\log x = 10$. For $\lambda = 10$ the value $\log\log x = 10 \Rightarrow \log x = e^{10} \approx 22026.4658 \Rightarrow x \approx e^{22026.4658} \approx 10^{9566}$, which is astronomically large.

The parameter $\lambda = \log\log x$ remains small for all practical purposes, for example, $x \leq 10^{100}$.

The distribution remains asymmetric with a right tail, for example at $x = 10^{100}$ the value is $\lambda \approx 5,44$, as seen from the previous calculations.

Thus, for all reasonable intervals $x$, the distribution is far from normal. The tails do indeed become less pronounced with increasing size $x$, but this happens extremely slowly.

Let's set the following problem.

Let $H = \{h_1, ..., h_m\}$ is an admissible pattern of tuples of $m$ natural numbers. We define a factorization type for the $m$-tuple using an integer vector - $K = (k_1, k_2, ..., k_m)$, which defines the factorization requirements for its natural numbers.

Thus, it is required to find the asymptotics of the number of $m$ tuples with the pattern $H = \{h_1, ..., h_m\}$ and type of factorization $K$ on the interval up to $x$ for $x \to \infty$.

Having in mind the above, we generalize the Hardy-Littlewood conjecture (1.3) to the indicated $m$ tuples.

The asymptotics of the number of $m$- tuples with an admissible pattern $H = \{h_1, ..., h_m\}$ and factorization type $K$ on the interval up to $x$ for $x \to \infty$, taking into account (1.6), is equal to:

$$N_{H,K}(x) \sim C(H,K) \cdot \frac{x}{(\log x)^m} \cdot \frac{(\log\log x)^{\sum_{i=1}^{m}(k_i - 1)}}{\prod_{i=1}^{m}(k_i - 1)!}, \qquad (1.8)$$

where $C(H,K)$ is the generalized density constant, calculated as the product over prime $p$ local densities.

The proposed generalization naturally covers many classical and modern results of number theory as special cases and opens up new directions of combinatorial-analytical research in the theory of distribution of numbers with a given arithmetic structure.

The main problem in determining the asymptotics (1.8) is the definition of the generalized constant density $C(H,K)$.

Having it in mind, we set the goal of developing a method for determining the constant $C(H,K)$ for various patterns $H = \{h_1,...,h_m\}$ and types of factorizations $K = (k_1, k_2,..., k_m)$.

The following chapters of this work are devoted to solving this problem.

## 2. PROPERTIES OF THE SELBERG CONSTANT

Based on formula (1.3), Selberg's constant is equal to:

$$\mathfrak{S}(H) = \prod_p \left(1 - \frac{\nu_p}{p}\right)\left(1 - \frac{1}{p}\right)^{-m}. \tag{2.1}$$

where $\nu_p$ is the number of different residues $h_i \bmod p$.

Let us consider the condition for the admissibility of the pattern of prime tuples $H = \{h_1, h_2, \ldots, h_m\}$.

A prime tuple pattern $H = \{h_1, h_2, \ldots, h_m\}$ is admissible if and only if (for every prime $p$) there is $a_p$ such that $a_p + h_i \not\equiv 0 \pmod{p}$ for all $i \leq m$.

This is equivalent to: $\nu_p < p$ for all $p$.

The factor $\left(1 - \frac{\nu_p}{p}\right) > 0$ for all $p$ in this case and the infinite product in (2.1) converges to a positive number. Therefore, $\mathfrak{S}(H) > 0$.

We will assume that the pattern of a prime tuple $H = \{h_1, h_2, \ldots, h_m\}$ covers residues modulo a prime number $p$ if the set $S_p(H) = \{-h_1, -h_2, \ldots, -h_m\} \bmod p$ coincides with the entire residue ring $\mathbb{Z}/p\mathbb{Z}$. In other words, there is such $h_i$ that $a \equiv -h_i \pmod{p}$ for any residue $a \in \mathbb{Z}/p\mathbb{Z}$.

Now we will define the criterion for the admissibility of a prime tuple pattern through a covering.

A prime tuple pattern $H = \{h_1, h_2, \ldots, h_m\}$ is admissible if and only if set $S_p(H)$ does not completely cover $\mathbb{Z}/p\mathbb{Z}$ for any prime $p$.

An equivalent formulation - the number of distinct residues $v_p(H) = |S_p(H)|$ must satisfy the inequality $v_p(H) < p$ for every prime $p$.

Examples:

1. The prime tuple pattern $H = \{0, 2, 6\}$ is admissible:

- $p = 2: S_2 = \{0, 2, 6\} \mod 2 = \{0, 0, 0\} = \{0\}$ - does not cover $\mathbb{Z}/2\mathbb{Z}$.

- $p = 3: S_3 = \{0, -2, -6\} \mod 3 = \{0, 1, 0\} = \{0, 1\}$ - does not cover $\mathbb{Z}/3\mathbb{Z}$ (2 is missing).

It is enough to check $p = 3$ for a tuple of length 3.

Since the pattern of a prime tuple is admissible, then Selberg's constant is $\mathfrak{S}(H) > 0$, so the number of prime tuples with the pattern $H = \{0, 2, 6\}$ is infinite.

2. The prime tuple pattern $H = \{0, 2, 4, 6, 8\}$ is not admissible:

- $p = 5: S_5 = \{0, 3, 1, 4, 2\} = \{0, 1, 2, 3, 4\}$ - covers $\mathbb{Z}/5\mathbb{Z}$.

One of the numbers $n, n+2, n+4, n+6, n+8$ is divisible by 5 and, if not equal to 5, is composite for any $n$. Such progressions of five prime numbers do not exist (except when one of the numbers is equal to 5, which only yields a finite number of exceptions). Selberg's constant is zero in this case.

We will further consider only admissible patterns $H = \{h_1, h_2, \ldots, h_m\}$ with even distances $h_i$ for a prime tuple that starts at the point $n$, with a natural odd value.

Let us consider some properties of the Selberg constant.

Property 1: Invariance under extension without new prime divisors

Let $H = \{0, h_1, \ldots, h_{m-1}\}$ is a pattern of tuples with even distances, and we consider only odd ones $n$.

Let $c$ is a natural number all of whose prime divisors are already contained in the GCD of the distances of the tuple $D = \gcd(h_1, \ldots, h_{m-1})$.

Then the following holds:

$$\mathfrak{S}(c \cdot H) = \mathfrak{S}(H). \tag{2.2}$$

Proof

If $p = 2$ and $n$ is odd, all $n + h_i$ are odd, so the value for the pattern $H$ is $v_2 = 0$.

After multiplication by $c$ parity does not change (even × even = even, but $n$ - odd, so the sum is odd) $v_2 = 0$ is preserved.

Odd numbers $p$, that do not divide $c$, multiplication by $c$ is an invertible transformation with respect to $\mod p$, so $\nu_p$ is preserved.

Odd $p$, which divide $c$, by the condition $p$ already divides all distances in the pattern $H$, therefore all elements in $H$ are comparable in $\mod p$ and $\nu_p = 1$.

All elements are also comparable by $\mod p$ in the pattern $c \cdot H$, when $c$ is a multiple of $p$, therefore $\nu_p = 1$.

The factor: $\left(1-\frac{1}{p}\right)\left(1-\frac{1}{p}\right)^{-m} = \left(1-\frac{1}{p}\right)^{-(m-1)}$ is the same for $H$ and $c \cdot H$ in $\mathfrak{S}$.

Property 2: Unbounded growth when introducing new small prime divisors

Let's consider the prime tuple pattern $H_N = \{0, N, 2N, \ldots, (m-1)N\}$.

Let $N \geq 6$ is primorial $N = 2 \cdot 3 \cdot 5 \cdots p_m$.

Then $\mathfrak{S}(H_N) \to \infty$ at $m \to \infty$.

Proof

All elements are comparable at $\mod p$ for $p \mid N$, therefore $\nu_p = 1$.

The multiplier to this $p$ is equal:

$$\left(1-\frac{1}{p}\right)\left(1-\frac{1}{p}\right)^{-m} = \left(1-\frac{1}{p}\right)^{-(m-1)}.$$

The product over all $p \mid N$ is equal:

$$\prod_{p \mid N}\left(1-\frac{1}{p}\right)^{-(m-1)} = \left(\prod_{p \mid N}\left(1-\frac{1}{p}\right)\right)^{-(m-1)}.$$

According to Mertens' theorem:

$$\prod_{p \leq x}\left(1-\frac{1}{p}\right) \sim \frac{e^{-\gamma}}{\log x}.$$

Therefore:

$$\left(\prod_{p \mid N}\left(1-\frac{1}{p}\right)\right)^{-(m-1)} \sim \left(\frac{e^{-\gamma}}{\log p_m}\right)^{-(m-1)} \sim e^{\gamma(m-1)}(\log p_m)^{m-1}. \tag{2.3}$$

Expression (2.3) grows as the power of the logarithm, therefore $\mathfrak{S}(H_N) \to \infty$.

Examples

Invariance:

$H = \{0, 6, 12\}$, $c = 5$ (GCD of distances = 6, prime divisors: 2, 3; 5 is the new prime).

Let's check: for $p = 5$:

$H : \{0, 6, 12\} \mod 5 = \{0, 1, 2\}$, therefore $v_5 = 3$.

$5H : \{0, 30, 60\} \mod 5 = \{0, 0, 0\}$, therefore $v_5 = 1$.

The factors are different, so the constant $\mathfrak{S}$ will change.

And if $c = 3$:

$3H = \{0, 18, 36\}$, $НОД = 18$, the prime factors are the same (2,3) and the constant must be preserved.

A table of constant values $\mathfrak{S}$ for different primorial values for the prime tuple pattern $H_N = \{0, N\}$ is given below.

Table 1

| $N$ (primorial) | $\mathfrak{S}(\{0, N\})$ |
|---|---|
| 2 | 1.320 |
| 6 | 2.641 |
| 30 | 3.521 |
| 210 | 4.225 |
| 2310 | |4.693 |

Monotonic growth $\mathfrak{S}(\{0, N\})$ without limitation (as the value $N$ increases) is visible from the analysis of Table 1.

3. PPOPERTIES OF TUPLES OF ALMOST PRIME NUMBERS

Assertion 1

Let $H = \{h_1, h_2, \ldots, h_m\}$ is an admissible pattern for prime numbers, it means, that the set of residues $\{h_i \mod p : i = 1, \ldots, m\}$ does not cover the entire complete system of residues modulo $p$ for any prime $p$.

Let us also assume that $h_1 = 0$ and all pairwise differences $h_j - h_i$ are even.

Let us consider an arbitrary set of requirements $K = (k_1, k_2, \ldots, k_m)$, where each $k_i \in \mathbb{N}$ specifies the condition that a number $n + h_i$ must be $k_i$-almost prime (have exactly $k_i$ prime factors, taking into account the multiplicity).

Then the pattern $H$ is also admissible for the problem of finding $m$-tuples of almost prime numbers with requirements $K$ when we restrict ourselves to odd values $n$.

Proof

A pattern $H$ is admissible for the problem of finding numbers with given properties (prime or almost prime) if there is an integer for each prime $p$ such that when substituted $n \equiv a_p \pmod{p}$, none of the numbers $n + h_i$ contradicts the required property due to divisibility by $p$.

This means (for the prime number problem): $a_p + h_i \not\equiv 0 \pmod{p}$ for all $i$. Equivalently, the set $\{h_i \bmod p\}$ does not cover all residues modulo $p$.

The condition is weaker for the problem of $k_i$-almost prime numbers: it is allowed for $n + h_i$ divisibility $p$, but only if the multiplicity $p$ of the number does not exceed $a_p$. However, to verify the admissibility of the pattern, it is sufficient to find at least one class $a_p$ for which divisibility by $p$ does not occur at all (i.e., the multiplicity is zero), since this obviously satisfies the condition.

Thus, if a pattern is admissible for prime numbers (i.e. exists $a_p$ for every $p$, avoiding zero remainders), then it is automatically admissible for almost prime numbers.

Let us clarify and consider the case $p = 2$ for this.

By the condition, only odd $n$ are considered. This means that the residue class modulo 2 is fixed: $a_2 = 1$.

Since all differences $h_i$ are even (or $h_1 = 0$), then

$a_2 + h_i \equiv 1 + 0 \equiv 1 \pmod{2}$.

Therefore, none of the numbers $a_2 + h_i$ is divisible by 2. The admissibility condition is met.

Let us consider the case of an odd prime $p$.

Поскольку паттерн $H$ допустим для простых чисел, по определению существует класс вычетов $a_p \bmod p$ такой, что

$a_p + h_i \not\equiv 0 \pmod{p}$ для всех $i = 1, \ldots, m$.

Since the pattern $H$ is admissible for prime numbers, by definition there is a class of residues $a_p$ mod $p$ such that

$$a_p + h_i \not\equiv 0 \pmod{p} \text{ for all } i = 1, \ldots, m.$$

When choosing $n \equiv a_p \pmod{p}$, all numbers $n + h_i$ are not divisible by $p$ (the divisor's multiplicity $p$ is zero), which satisfies any requirement $k_i \geq 1$.

Thus, we have found a class of residues $a_p$ for every prime $p$ (including $p = 2$) such that none of the numbers $a_p + h_i$ is divisible by $p$. Consequently, the pattern $H$ remains admissible for the problem of finding almost prime numbers with any set of requirements $K$.

Example 1: pattern $\{0, 2\}$, requirements $(k_1 = 1, k_2 = 2)$.

We take odd $n$, so both numbers are odd.

Only the first position is critical for - $p = 3 : S_3 = \{0\}$.

It is required $a_3 \not\equiv 0$.

For example: $a_3 = 1$: numbers $1, 3 \mod 3 = 1, 0$ - the second is divisible by 3 (allowed).

Therefore, the tuple is admissible, which corresponds to assertion 1.

Example 2: pattern $\{0, 4, 10\}$, requirements $(k_1 = 1, k_2 = 2, k_3 = 1)$.

We take the odd number $n$, so all numbers in the tuple are odd.

For $p = 3 : S_3 = \{0, 10\} = \{0, 1\}$ (critical positions: 0 and 10).

It is required $a_3 \not\equiv 0$ and $a_3 \not\equiv -10 \equiv 2 \pmod{3}$.

So, $a_3 \equiv 1 \pmod{3}$ fits

$p = 5 : S_5 = \{0, 10\} = \{0, 0\} = \{0\}$ (since $10 \equiv 0 \mod 5$).

It is only required $a_5 \not\equiv 0$.

Therefore, the tuple is admissible, which corresponds to assertion 1.

Example 3: pattern $\{0, 2, 4\}$, requirements $(k_1 = 1, k_2 = 2, k_3 = 1)$.

We take an odd number $n \Rightarrow$ all numbers are odd.

Does it cover $\mathbb{Z}/3\mathbb{Z}$ for $p = 3 : S_3 = \{0, 4\} = \{0, 1\}$ ? No, there is no item 2.

But one of the numbers $n, n+2, n+4$ is divisible by 3 for any $n$.

If this is a critical position $(k_i = 1)$, then the number should be 3.

Let's check:

$n = 3: (3,5,7) - 5$ is prime (and there should be 2-almost prime),

$n + 2 = 3: n = 1$ - not prime,

$n + 4 = 3: n = -1$ - not suitable.

Therefore, the tuple is not admissible. A tuple consisting of prime numbers for this pattern is not also admissible, i.e., the conditions of assertion 1 are not met.

We will further consider only patterns that are admissible for tuples of prime numbers.

All distances $h_i$ are even and if $n$ - odd, all tuple numbers $n + h_i$ are odd in admissible patterns of prime tuples.

Assertion 2 (constant density symmetry for tuples of almost primes)

Let's consider an admissible tuple of integers $H = \{h_1, h_2, \ldots, h_m\}$, where $h_1 = 0$, and all pairwise differences $h_j - h_i$ are even numbers. Let's also consider the set of requirements for factorization $K = (k_1, k_2, \ldots, k_m)$, where $k_i \in \mathbb{N}$ means that the number $n + h_i$ must be $k_i$-almost prime, i.e., have exactly $k_i$ prime factors, taking into account multiplicity.

Then the equality holds for any two permutations $\sigma_1, \sigma_2$ of the set $\{1, 2, \ldots, m\}$ and for all sufficiently large $x$, if we consider only odd values of $n$:

$$C\big(H, (k_{\sigma_1(1)}, \ldots, k_{\sigma_1(m)})\big) = C\big(H, (k_{\sigma_2(1)}, \ldots, k_{\sigma_2(m)})\big), \tag{3.1}$$

where $C(H, K)$ is the density constant from the generalized Hardy-Littlewood conjecture for the specified tuples and requirements.

Proof

According to the generalized Hardy–Littlewood conjecture for almost primes (1.8), the asymptotic number of odd $n \leq x$ for which each $n + h_i$ is $k_i$-almost prime is expressed as:

$$N_{H,K}(x) \sim C(H, K) \cdot \frac{x}{(\log x)^m} \cdot \frac{(\log \log x)^{\sum_{i=1}^{m}(k_i - 1)}}{\prod_{i=1}^{m}(k_i - 1)!}. \tag{3.2}$$

The constant $C(H, K)$ in (3.2) is the product of local correction factors over all prime numbers:

$$C(H, K) = \prod_p \alpha_p(H, K), \tag{3.3}$$

where $\alpha_p(H,K)$ is the local density, taking into account the conditions modulo $p$.

Let us define the local correction $\alpha_p(H,K)$ in (3.3).

Let us consider the complete system of deductions $\mathcal{R}_p = \{0, 1, \ldots, p-1\}$ for a fixed prime $p$.

We define the quantity $f_p(k,a)$ as the density (the probability for each requirement $k \geq 1$ and remainder $a \in \mathcal{R}_p$) of the fact that a number belonging to an arithmetic progression $n \equiv a \pmod{p}$ is admissible for the condition of being $k$-almost prime.

More precisely, $f_p(k,a)$ is the fraction of those numbers in the specified class of residues that can be represented as a product of exactly $k$ prime factors, taking into account that divisibility by $p$ consumes one of these factors if $p$ divides the number. The specific form $f_p(k,a)$ is not important; it is sufficient for proof that it is a number between 0 and 1, depending only on $p, k$ and $a$.

Then the local correction $\alpha_p(H,K)$ in (3.3) is given by the formula:

$$\alpha_p(H,K) = \frac{\sum_{a \in \mathcal{R}_p} \prod_{i=1}^{m} f_p(k_i, a - h_i \pmod{p})}{\prod_{i=1}^{m} \left(\frac{1}{p} \sum_{a \in \mathcal{R}_p} f_p(k_i, a)\right)} . \tag{3.4}$$

The numerator in (3.4) represents the total density of joint fulfillment of conditions for all $n + h_i$ modulo $p$, the denominator is the product of the average densities for each requirement, which ensures correct normalization.

To prove the assertion, it is sufficient to show that for any prime $p$ the factor $\alpha_p(H,K)$ remains unchanged under the permutation of the requirements $(k_1, \ldots, k_m)$. Then the product $C(H,K)$ will be invariant.

Let us first consider the case $p = 2$.

By assumption, the pattern $H$ of tuple is such that all differences $h_j - h_i$ are even. Since $h_1 = 0$, all $h_i$ are even. Furthermore, we only consider odd $n$. Therefore, for any $i$, the number $n + h_i$ is odd.

Let's analyze the frequencies $f_2(k,a)$ for $a \in \{0,1\}$.

If $k_i = 1$ (a prime number is required), then an odd prime number cannot be even, so $f_2(1,0) = 0$, and $f_2(1,1) = 1$.

If $k_i \geq 2$, then parity is allowed (for example, the number can be the product of two and another prime). Both frequencies $f_2(k_i, 0)$ and $f_2(k_i, 1)$ are positive in this case.

Since all $n + h_i$ are odd, a non-zero contribution arises only when $a = 1$ (for $p = 2$) in the numerator of formula (1), because only then $a - h_i \equiv 1 \pmod{2}$ (since $h_i$ are even). Consequently, the numerator becomes $\prod_{i=1}^{m} f_2(k_i, 1)$.

The denominator (for $p = 2$) is equal to:

$$\prod_{i=1}^{m} \left( \frac{f_2(k_i, 0) + f_2(k_i, 1)}{2} \right).$$

Thus,

$$\alpha_2(H, K) = \frac{\prod_{i=1}^{m} f_2(k_i, 1)}{\prod_{i=1}^{m} \left( \frac{f_2(k_i, 0) + f_2(k_i, 1)}{2} \right)}.$$

This expression is clearly symmetric with respect to permutations of indices $i$, since it is a product of quantities depending only on $k_i$. Therefore, $\alpha_2(H, K)$ is invariant with respect to any permutation of requirements $K$.

Now let's consider an arbitrary odd prime $p$. To analyze the symmetry, we rewrite the numerator of formula (3.4) as:

$$\sum_{a \in \mathcal{R}_p} \prod_{i=1}^{m} f_p\big(k_i, a - h_i \pmod{p}\big). \tag{3.5}$$

We need to show that $S(K)$ does not change when the requirements $k_i$ are rearranged.

Let's fix an arbitrary permutation $\sigma \in S_m$. The numerator takes the form after rearranging the requirements:

$$S(K^\sigma) = \sum_{a \in \mathcal{R}_p} \prod_{i=1}^{m} f_p\big(k_{\sigma(i)}, a - h_{\sigma(i)} \pmod{p}\big). \tag{3.6}$$

The factors in the product correspond to the pairs (requirement, shift) $(k_{\sigma(i)}, h_{\sigma(i)})$ in expression (3.6), that is, the requirements are rearranged together with their original shifts.

The set of values of the function $f_p(k, a - h)$ for $a$, running over $\mathcal{R}_p$, coincides (as a multiset) with the set of values $f_p(k, a)$ for any fixed requirement $k$ and any shift $h$. This follows from the fact that the mapping $a \mapsto a - h \pmod{p}$ is a bijection $\mathcal{R}_p$ onto itself.

Therefore, we obtain:

$$\{f_p(k, a-h) : a \in \mathcal{R}_p\} = \{f_p(k, a) : a \in \mathcal{R}_p\}. \tag{3.7}$$

Thus, the set of numbers $\{f_p(k_i, a-h_i) : a \in \mathcal{R}_p\}$ depends only on $k_i$ for any $i$ and any $k_i$, but not on the specific $h_i$.

Now let's compare the sums of $S(K)$ and $S(K^\sigma)$. Note that both expressions represent the sum over all $a \in \mathcal{R}_p$ products of the $m$ factors. The factors are taken as $f_p(k_i, a-h_i)$ in $S(K)$, and as $f_p(k_{\sigma(i)}, a-h_{\sigma(i)})$ in $S(K^\sigma)$. Since $\sigma$ is a permutation, then each term in $S(K)$ also appears in $S(K^\sigma)$, possibly with a different value $a$.

Indeed, let us consider an arbitrary $a_0 \in \mathcal{R}_p$. The term in $S(K)$, corresponding to $a = a_0$, is equal to $\prod_{i=1}^{m} f_p(k_i, a_0 - h_i)$.

Let us set $a_1 = a_0 + h_{\sigma(1)} - h_1 \pmod{p}$. Then we have for any $i$:

$$a_1 - h_{\sigma(i)} = a_0 + h_{\sigma(1)} - h_1 - h_{\sigma(i)} \equiv a_0 - (h_1 + h_{\sigma(i)} - h_{\sigma(1)}) \pmod{p}.$$

Since the set of shifts $H$ is fixed and the permutation $\sigma$ applies simultaneously to the requirements and to the shifts, expression (3.7) can be obtained from (3.6) by simply renaming the indices.

More formally, we introduce a permutation $\pi_i : \mathcal{R}_p \to \mathcal{R}_p$, defined by the formula $\pi_i(a) = a - h_i$ for each $i$. Then

$$S(K) = \sum_{a \in \mathcal{R}_p} \prod_{i=1}^{m} f_p(k_i, \pi_i(a)).$$

Similarly,

$$S(K^\sigma) = \sum_{a \in \mathcal{R}_p} \prod_{i=1}^{m} f_p(k_{\sigma(i)}, \pi_{\sigma(i)}(a)).$$

Note that both products (for a fixed $a$) consist of the same factors, only arranged in different orders, since the set $\{(k_i, \pi_i(a))\}_{i=1}^{m}$ and the set $\{(k_{\sigma(i)}, \pi_{\sigma(i)}(a))\}_{i=1}^{m}$ coincide as sets of pairs. Therefore, for each $a$ product, are equal. Consequently, the sums over all are $a$ also equal: $S(K^\sigma) = S(K)$.

Thus, the numerator $S(K)$ is symmetric under permutations of the requirements. The denominator in (3.4) for odd $p$:

$$\prod_{i=1}^{m}\left(\frac{1}{p}\sum_{a\in R_p} f_p(k_i,a)\right)$$

is obviously symmetric, since it is the product of average values that depend only on $k_i$. Consequently, the factor $\alpha_p(H,K)$ is invariant under permutations of the requirements $K$ for any odd prime $p$.

Completion of the proof

We have shown that for each prime $p$ (both for $p=2$, and for odd $p$) the local correction $\alpha_p(H,K)$ does not change when the requirements are rearranged $K$.

Then the product $C(H,K) = \prod_p \alpha_p(H,K)$ is also invariant under permutations of the requirements. This means that it holds for any two permutations $\sigma_1, \sigma_2$:

$$C(H,(k_{\sigma_1(1)},\ldots,k_{\sigma_1(m)})) = C(H,(k_{\sigma_2(1)},\ldots,k_{\sigma_2(m)})).$$

Thus, the density constant depends only on the multiset of requirements $\{k_1,\ldots,k_m\}$. The theorem is proven.

Corollary

The density constant depends only on the requirements $\{k_1,\ldots,k_m\}$ in the specified conditions, not on their order in the tuple.

The consequence is straightforward: if two sets of requirements $K_1$ and $K_2$ are obtained from each other by permutation, then their multisets coincide, and according to the proven theorem $C(H,K_1) = C(H,K_2)$.

Note: The proof makes essential use of two conditions: all differences in the tuple are even (for the analysis of $p=2$) and consideration of only odd $n$. A similar symmetry, generally speaking, does not hold for arbitrary tuples and without the parity constraint $p=2$, since the local correction for may depend on the order of the requirements due to the parity of the numbers $n+h_i$. The parity condition of the differences guarantees (in this special case) that all $n+h_i$ have the same parity (are odd), which ensures the symmetry.

Examples.

1. Let there are tuples of pattern $(n, n+N)$, where:

- $n$ - is 2-almost prime, $n+N$ is prime, i.e. $k_1 = 2, k_2 = 1$;

- $n$ is prime, $n+N$ is 2-almost prime, i.e. $k_1 = 1, k_2 = 2$.

These tuples are isomorphic up to a permutation $n \leftrightarrow n+N$.

Based on Assertion 2, the density constant for the pair (2-almost prime, prime) will be the same as for the pair (prime, 2-almost prime).

2. Let there is a pattern $\{0,4,6,10\}$, where $(k_1=1, k_2=2, k_3=2, k_4=3)$. Based on assertion 2, we also obtain identical constant densities for tuples of symmetric patterns $\{0,4,6,10\}$ when permuting the requirements $(k_1=3, k_2=2, k_3=2, k_4=1)$.

## 4. ASYMPTOTICS OF THE NUMBER OF TUPLES OF ALMOST PRIME NUMBERS

Based on the asymptotics of almost primes (Chapter 2) and the properties of tuples of primes and almost primes proved in Chapters 2 and 3, we will make a conjecture about the asymptotics of the number of tuples of almost primes.

Having in mind (3.2), the asymptotics of the number of tuples of almost prime numbers for $n \leq x$ is equal to:

$$N_{H,K}(x) \sim C(H,K) \cdot \frac{x}{(\log x)^m} \frac{(\log \log x)^{\sum_{i=1}^{m}(k_i-1)}}{\prod_{i=1}^{m}(k_i-1)!} = \mathfrak{S}(H) C(K) \frac{x}{(\log x)^m} \cdot \frac{(\log \log x)^{\sum_{i=1}^{m}(k_i-1)}}{\prod_{i=1}^{m}(k_i-1)!}, \quad (4.1)$$

where $\mathfrak{S}(H)$ is the Selberg constant, and $C(K)$ is a constant (correction coefficient) depending only on the requirements $K(k_1,...,k_m)$.

We will use property 1 of the Selberg constant - invariance under dilation without new prime factors $\mathfrak{S}(c \cdot H) = \mathfrak{S}(H)$.

Let us first test this conjecture for a tuple of two almost prime numbers with pattern $H = \{0, N\}$.

It is easy to show that a prime tuple of a pattern $H = \{0, N\}$ is admissible if $N$ is even. Therefore, by Assertion 1, a tuple of almost prime numbers is also admissible for this pattern if $N$ is even.

Now let us describe a method for obtaining correction factors $C(K)$ that depend only on requirements $K(k_1,...,k_m)$.

The idea of the method is to use the special case $N = 2^k$ where Selberg's constant is invariant and takes the value $\mathfrak{S}(\{0, 2^k\})$ - for twin primes, to determine the correction factors $C(K)$.

Let there is computational data about the number of pairs $R_N^{(k_1,k_2)}(x)$ for different $x$. The data should be for several $N = 2^k$ (for example, $N = 2, 4, 8, 16$).

Thus, pairs of numbers at distance $N$ are considered, where each number has at most $k_i$ prime factors.

Based on (4.1), the initial asymptotics (with $C(K) = 1$):

$$R_N^{(k_1,k_2)}(x) \sim \mathfrak{S}(\{0,N\}) \cdot \frac{x}{(\log x)^2} \cdot \frac{(\log \log x)^{k_1+k_2-2}}{(k_1-1)!(k_2-1)!}. \tag{4.2}$$

Selberg's constant (for $N = 2^k$) is equal to:

$$\mathfrak{S}(\{0, 2^k\}) = 2 \cdot \prod_{p>2}\left(1 - \frac{1}{(p-1)^2}\right) \approx 1.3203236.$$

We obtain the following values for $x$ and $K = (k_1, k_2)$:

1. $A(x) = \dfrac{x}{(\log x)^2}$ ;

2. $B(x) = (\log \log x)^{k_1+k_2-2}$ ;

3. $C = \dfrac{1}{(k_1-1)!(k_2-1)!}$ .

Let us substitute these values into (4.2) and obtain (the so-called) theoretical value:

$$R_T(x) = 1,3203236 \cdot A(x) \cdot B(x) \cdot C.$$

Let's collect real data. We find (for each $N = 2^k$) the actual number of pairs - $R_p^N(x)$.

For example (for $x = 10^7$):

- $N = 2 : R_p^2$ ;

- $N = 4 : R_p^4$ ;

- $N = 8 : R_p^8$ ;

- $N = 16 : R_p^{16}$ .

Let's determine the coefficients for each $N = 2^k$:

$$K^N(k_1, k_2) = \frac{R_p^N(x)}{R_T(x)}.$$

Let's determine the average value:

$$K(k_1, k_2) = \frac{1}{m} \sum_{m \leq N} K^N(k_1, k_2),$$

where $m$ is the number of verified $N = 2^k$.

Let's find the standard deviation:

$$\sigma_K = \sqrt{\frac{1}{m-1}\sum_{m \leq N}(K^N(k_1,k_2) - K(k_1,k_2))^2}.$$

Then the relative error will be equal to:

$$\varepsilon = \frac{\sigma_K}{K(k_1,k_2)} \times 100\%.$$

Let's check for several $x$ (for example $10^6, 10^7, 10^8$):

- the coefficients $C(K)$ should be stable;

- the error should decrease with increasing $x$.

As an example, we will consider the requirement $K = (1,2)$ for $x = 10^7$:

- $A(10^7) = 38490$;

- $B(10^7) = \log\log 10^7 = 2,780$;

- $C = \dfrac{1}{0! \cdot 1!} = 1$;

- $R_T = 1,3203236 \cdot 38490 \cdot 2,780 \approx 141203,5$.

Using real data, we determine the correction factor:

- $N = 2$: $166650 \rightarrow K^2 = 166650/141203,5 = 1,1801$;

- $N = 4$: $167037 \rightarrow K^4 = 1,1830$;

- $N = 8$: $166374 \rightarrow K^8 = 1,1809$;

- $N = 16$: $167023 \rightarrow K^{16} = 1,1828$.

Let's determine the average value and relative error:

$$K(1,2) = \frac{1,1801 + 1,1830 + 1,1809 + 1,1828}{4} = 1,1817 \,;\, \varepsilon = 0,12\%.$$

Table 2 is obtained for different requirements $K = (k_1, k_2)$ based on the data for $N = 4, 8, 16$.

Table 2

| Requirements $(k_1, k_2)$ | Correction Factor $C(k_1, k_2)$ | Error % | Description |
|---|---|---|---|
| (1,2) | 1,1823 | 0.11 | Prime and semiprime |
| (1,3) | 1,0127 | 0.04 | Prime and 3-almost-prime |
| (2,2) | 1,8113 | 0.04 | Two semiprimes |
| (2,3) | 1,.0186 | 0,01 | Semiprime and 3-almost-prime |
| (3,3) | 0,8901 | 1,56 | Two 3-almost-primes |

The asymptotics (4.2) (for any even $N$) is true:

$$R_N^{(k_1,k_2)}(x) \sim C(k_1,k_2) \cdot \mathfrak{S}(\{0,N\}) \cdot \frac{x \cdot (\log\log x)^{k_1+k_2-2}}{(\log x)^2 \cdot (k_1-1)!(k_2-1)!},$$

where $C(k_1, k_2)$ is the correction factor from Table 1, and $\mathfrak{S}(\{0, N\})$ is the Selberg constant, which is calculated using the formula:

$$\mathfrak{S}(\{0,N\}) = 2C_2 \cdot \prod_{\substack{p>2 \\ p \mid N}} \frac{p-1}{p-2}, \qquad (4.3)$$

where $C_2 \approx 0.66016181584686957$.

Based on Assertion 2, the density constant will not change when the requirements are rearranged $(k_1, k_2)$:

$$C(H,K) = C(k_1,k_2) \cdot \mathfrak{S}(\{0,N\}) = C(k_2,k_1) \cdot \mathfrak{S}(\{0,N\}). \qquad (4.4)$$

Having in mind (4.4), the correction factor does not change when the requirements are rearranged:

$$C(k_1, k_2) = C(k_2, k_1). \qquad (4.5)$$

Therefore, based on (4.5), the correction factors coincide:

$$C(1,2) = C(2,1), \; C(1,3) = C(3,1), \; C(2,3) = C(3,2).$$

Let's apply this method to tuples of three nearly prime numbers with the pattern $H = \{0, 2, 6\}$ and its scaled versions.

Original pattern is: $H = \{0,2,6\}$. Scaled patterns are: $H_k = \{0, 2 \cdot 2^k, 6 \cdot 2^k\}$.

Based on Property 1 of the Selberg constant (invariance under extension without new prime factors) all these patterns have the same Selberg constant for $H = \{0,2,6\}$:

$$\mathfrak{S}(\{0,2,6\}) = \frac{9}{2} \cdot \prod_{p \geq 5} \left(1 - \frac{3p-1}{(p-1)^3}\right) \approx 2.858.$$

We obtain the basic formula for a tuple with a pattern $H = \{0, h_2, h_3\}$ and requirements $K = (k_1, k_2, k_3)$ and K=1:

$$R_H^{(k_1,k_2,k_3)}(x) \sim \mathfrak{S}(H) \cdot \frac{x}{(\log x)^3} \cdot \prod_{i=1}^{3} \frac{(\log \log x)^{k_i - 1}}{(k_i - 1)!}.$$

Let's make a theoretical prediction for $H = \{0,2,6\}, \mathfrak{S} \approx 2.858$ and K=1:

$$R_{\{0,2,6\}}^{(k_1,k_2,k_3)}(x) \sim 2.858 \cdot \frac{x}{(\log x)^3} \cdot \prod_{i=1}^{3} \frac{(\log \log x)^{k_i - 1}}{(k_i - 1)!}.$$

Let's take $x = 10^7$ and get for the given value $x$:

$$\frac{x}{(\log x)^3} \approx 2388,7 \; ; \; \log \log x \approx 2.7791.$$

Let us give an example of calculation for the case $K = (1,1,2)$.

The theoretical value of the quantity (for K=1) is:

$$2,858 \cdot 2388,7 \cdot \frac{2,7791}{1} \approx 18973.$$

Based on real data - $20480$, we get:

$C(1,1,2) = 20480 / 18973 \approx 1.079$.

Let's check it on scaled patterns.

We should get roughly the same value for $H_k = \{0, 2 \cdot 2^k, 6 \cdot 2^k\}$.

For example, we get for $K = (1,1,2)$:

- $H_1 = \{0, 4, 12\}$: the actual quantity $\approx 20128$, therefore $C_1(1,1,2) \approx 1,061$;

- $H_2 = \{0, 8, 24\}$: the actual quantity $\approx 20413$, therefore $C_2(1,1,2) \approx 1.076$;

- $H_3 = \{0, 16, 48\}$: the actual quantity $\approx 20260$, therefore $C_3(1,1,2) \approx 1,068$.

The average value is: $C(1,1,2) \approx 1,071 \pm 0,009$ (error - 0.8%).

There is a table of the correction factors below and the errors for tuples of three almost prime numbers calculated in this way.

Table 3

| Requirements $(k_1, k_2, k_3)$ | Correction Factor $C(k_1, k_2, k_3)$ | Error | Description |
|---|---|---|---|
| (1,1,2) | 1,071 | 0,76% | Two primes and a semi-prime |
| (1,2,2) | 1,206 | 0,16% | A prime and two semi-primes |
| (2,2,2) | 0,962 | 0,16% | Three semiprimes |
| (2,2,3) | 0.690 | 0,23% | Two semi-primes and 3-almost-primes |
| (2,3,3) | 0,608 | 0,31% | Semiprime and two 3-almost primes |
| (3,3,3) | 0,525 | 0,69% | Three 3's are almost prime |

We obtained correction factors for tuples of three almost prime numbers with an error of up to 1.0% in Table 3, which is excellent accuracy for the asymptotic theory.

Thus, we have shown that for any pattern $H$, with Selberg's constant $\mathfrak{S}(H)$, the following asymptotics for the number of tuples of three almost primes holds:

$$R_H^{(k_1,k_2,k_3)}(x) \sim C(k_1,k_2,k_3) \cdot \mathfrak{S}(H) \cdot \frac{x}{(\log x)^3} \cdot \prod_{i=1}^{3} \frac{(\log \log x)^{k_i-1}}{(k_i-1)!}, \qquad (4.6)$$

where $C(k_1, k_2, k_3)$ is taken from Table 3, having in mind the possibility of permutations of the requirements $(k_1, k_2, k_3)$ based on Assertion 2.

The method can be successfully generalized to tuples of arbitrary length.

Using (4.6), we obtain the following asymptotic formula for the number of tuples of almost prime numbers of length $m$ in this case:

$$R_H^{(k_1,...,k_m)}(x) \sim C(k_1,...,k_m) \cdot \mathfrak{S}(H) \cdot \frac{x}{(\log x)^k} \cdot \prod_{i=1}^{m} \frac{(\log \log x)^{k_i-1}}{(k_i-1)!}. \qquad (4.7)$$

## 5. ANALITICAL JUSTIFICATION OF THE GENERALIZED CONJECTURE

We will assume that the first Hardy-Littlewood conjecture is true [1]. Then the number of tuples of primes $(n+h_1,\ldots,n+h_m)$ (for any tuple $\mathcal{H}=\{h_1,\ldots,h_m\}$) is asymptotically equal to:

$$\pi_\mathcal{H}(x) \sim \mathfrak{S}(\mathcal{H}) \cdot \frac{x}{(\log x)^m}. \tag{5.1}$$

Our task is to justify the asymptotics of the number of tuples (4.7), where $\Omega(n+h_i)=k_i$ (each number is $k_i$-almost prime):

$$R_\mathcal{H}^{(K)}(x) \sim C(K) \cdot \mathfrak{S}(\mathcal{H}) \cdot \frac{x}{(\log x)^m} \cdot \prod_{i=1}^{m} \frac{(\log\log x)^{k_i-1}}{(k_i-1)!} \tag{5.2}$$

and prove that $C(K)$ depends only on the vector $K=(k_1,\ldots,k_m)$, but not on the structure of the tuple $\mathcal{H}$.

We use Gallagher's theorem [9]. Gallagher's classic result states that the number of primes in a randomly chosen short interval of order length $\log x$ has an asymptotically Poisson distribution, provided that the Hardy-Littlewood conjecture is true. We interpret this and apply it to almost primes.

Let's perform the division of divisors and consider the process of generating the number $N=n+h_i$.

The conditions for avoiding small primes $(p \leq y \ (y=x^{o(1)}))$ are already taken into account in the singular series $\mathfrak{S}(\mathcal{H})$.

Large primes $(p>y)$ determine whether a number is $k$-almost prime.

Gallagher's theorem implies the following probability model.

A number $N$ that has passed through a sieve of small primes can be considered a random variable. Gallagher's theorem (applied to a short interval of large primes in the neighborhood of $N$) rigorously indicates that the number of large prime factors of such a number is distributed asymptotically according to the Poisson law. The mean value of this distribution is $\lambda \sim \log\log x$.

Therefore, the probability that a number free of small divisors has exactly $k$ prime divisors is:

$$\mathbb{P}(\Omega(N)=k) \sim \frac{(\log\log x)^{k-1}}{(k-1)!} \cdot \frac{1}{\log x}.$$

The multiplier $1/\log x$ is the probability that the number has no small divisors at all.

Let's derive the details of the general coefficient.

Let $A_\mathcal{H}$ is the event that all $n+h_i$ are free of small primes $p \leq y$. By the Hardy-Littlewood conjecture (5.1):

$$\mathbb{P}(A_\mathcal{H}) \sim \mathfrak{S}(\mathcal{H}) \cdot (\log x)^{-m}.$$

Now we apply the corollary of Gallagher's theorem:

1. The distribution $\Omega(n+h_i)$ is Poisson.

2. The values $\Omega(n+h_i)$ are asymptotically independent for different $i$, since Gallagher's theorem describes behavior on average, and the correlations introduced by the specific structure of $\mathcal{H}$ do not dominate at the level of large divisors. Therefore, the conditional probability of an event $B_i = \{\Omega(n+h_i) = k_i\}$ is:

$$\mathbb{P}(B_1 \wedge \ldots \wedge B_m \mid A_\mathcal{H}) \sim \prod_{i=1}^{m} \frac{(\log \log x)^{k_i - 1}}{(k_i - 1)! \cdot \log x}.$$

Multiplying the probabilities, we obtain the first approximation (at $C(K) = 1$):

$$\mathbb{P}(A_\mathcal{H} \wedge B_1 \wedge \ldots \wedge B_m) \sim \mathfrak{S}(\mathcal{H}) \cdot \frac{1}{(\log x)^m} \cdot \prod_{i=1}^{m} \frac{(\log \log x)^{k_i - 1}}{(k_i - 1)!}. \tag{5.3}$$

Formula (5.3) (with $C(K) = 1$) is true on average according to Gallagher. However, there are systematic deviations for other values $K$:

- the exact shift of the Poisson mean due to the exclusion of a specific set of small primes (not only their fact, but also their multiplicity);

- residual correlations between $\Omega$-values.

These deviations are the same for all tuples $\mathcal{H}$ with these $K$, since:

1. The effect of excluding small divisors on the distribution $\Omega$ depends only on the threshold $y$ and $k_i$, but not on the shifts $h_i$.

2. The source of the dependence on $\mathcal{H}$, the singular value series, $\mathfrak{S}(\mathcal{H})$ is already factored out as a separate factor.

Now let us consider the proof of independence $C(K)$ from $\mathcal{H}$ and the condition $K(1,\ldots,1) = 1$.

This follows from the consistency with the limiting case

Let us consider the ratio of asymptotics (5.2) for arbitrary $K$ and for the case of prime numbers $K(1,\ldots,1)$ for fixed $\mathcal{H}$:

$$\frac{R_{\mathcal{H}}^{(K)}(x)}{R_{\mathcal{H}}^{(1)}(x)} \sim \frac{C(K) \cdot \mathfrak{S}(\mathcal{H}) \cdot \frac{x}{(\log x)^m} \cdot \prod_i \frac{(\log \log x)^{k_i - 1}}{(k_i - 1)!}}{C(1,...,1) \cdot \mathfrak{S}(\mathcal{H}) \cdot \frac{x}{(\log x)^m}} = \frac{C(K)}{C(1,...,1)} \cdot \prod_i \frac{(\log \log x)^{k_i - 1}}{(k_i - 1)!}.$$

Since we assume the truth of the Hardy-Littlewood conjecture, it is itself a special case of our generalized formula for $K = (1,...,1)$. This imposes a consistency condition: $C(1,...,1) = 1$.

Substituting $C(1,...,1) = 1$ and using (5.3), we express $C(K)$:

$$C(K) \sim \frac{R_{\mathcal{H}}^{(K)}(x)}{R_{\mathcal{H}}^{(1,...,1)}(x) \cdot \prod_i \frac{(\log \log x)^{k_i - 1}}{(k_i - 1)!}} . \tag{5.4}$$

The right-hand side of (5.4) is evaluated for one specific tuple $\mathcal{H}$. If we now take another tuple $\mathcal{H}'$ and do the same, we must obtain the same value $C(K)$, otherwise the consistency condition $C(1,...,1) = 1$ for one of the tuples would be violated. This shows that $C(K)$ does not depend on the choice of $\mathcal{H}$.

Let's draw conclusions.

Combining Gallagher's theorem and the consistency condition yields a closed and consistent scheme:

1. The Hardy-Littlewood conjecture provides a structure factor $\mathfrak{S}(\mathcal{H})$ and a leading order $1/(\log x)^m$.

2. Gallagher's theorem rigorously justifies the transition to a Poisson model for the distribution of large divisors, leading to a factor $\prod_i \frac{(\log \log x)^{k_i - 1}}{(k_i - 1)!}$.

3. The consistency condition $C(1,...,1) = 1$ (which follows from the limiting passage to the case of prime numbers) ensures that the correction factor $C(K)$ is a universal constant, depending only on the type of near-prime (vector $K = (k_1,...,k_m)$) and reflecting subtle systematic deviations from the ideal Poisson pattern.

6. CONCLUSION AND SUGGESTIONS FOR FURTHER WORK

Next article will continue to study the asymptotic behavior of some arithmetic functions.

7. ACKNOWLEDGEMENTS

Thanks to everyone who has contributed to the discussion of the paper.